\documentclass[a4paper,12pt]{article}
\usepackage{graphicx,amssymb,color}
\usepackage{epsfig}
\usepackage{color}

\begin{document}

\centerline{\LARGE \bf Travelling waves for a bistable reaction-diffusion}

\medskip

\centerline{\LARGE \bf equation with delay}

\vspace*{1cm}

\centerline{\bf S. Trofimchuk$^1$, V. Volpert$^{2,3,4}$}
%
\bigskip

\centerline{$^1$
Instituto de Matem´atica y Fisica, Universidad de Talca, Casilla 747, Talca, Chile}


\centerline{$^2$
Institut Camille Jordan, UMR 5208 CNRS, University Lyon 1, 69622 Villeurbanne, France}

\centerline{$^3$
INRIA
Team Dracula, INRIA Lyon La Doua, 69603 Villeurbanne, France}

\centerline{$^4$
Peoples' Friendship University of Russia, ul. Miklukho-Maklaya 6, Moscow, 117198 Russia}

%

\vspace*{1cm}

\noindent
{\bf Abstract.}
The paper is devoted to a reaction-diffusion equation with delay arising in modelling the immune response.
We prove the existence of travelling waves in the bistable case using the Leray-Schauder method.  In difference with the previous works, we do not assume here quasi-monotonicity of the delayed reaction term.

\vspace*{0.2cm}

\noindent
{\bf Key words:} reaction-diffusion equation, time delay, bistable case, wave existence

\vspace*{4mm}


\section{Introduction}

\begin{figure}[ht!]
\centerline{\includegraphics[scale=0.5]{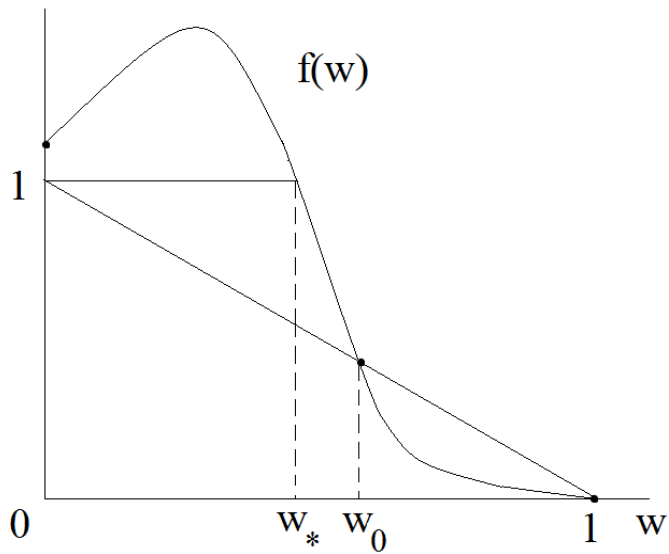}}
\caption{The typical form of the function $f(v)$. It is growing for small $v$ and decreasing for large $v$.}
\label{f7}
\end{figure}

In this work we study the existence of travelling waves for the reaction-diffusion equation with delay:
\begin{equation}\label{m0}
  \frac{\partial v}{\partial t} = D \frac{\partial^2 v}{\partial x^2} + k v (1-v) - f(v_\tau) v .
\end{equation}
Here $v=v(x,t)$, $v_\tau=v(x,t-\tau)$, the function $f(v_\tau)$ will be specified below.
This equation models the spreading of viral infection in tissues such as spleen or lymph nodes (see  \cite{PLOS}).
The first term in the right-hand side of this equation describes virus diffusion, the second term its production
and the last term its elimination by the immune cells.  The parameter $D$ is the diffusion coefficient (or diffusivity) and $k$ stands for the replication rate constant.  In  the sequel, without loss of  generality, we can assume that $D=k=1 $. The  parameterised function $f(v_\tau)$  (where $v_\tau$ is the concentration  of virus some time $\tau$ before) characterises the virus induced clonal expansion of T cells, i.e. the number  and function of these cells  upon their maturation during some time  $\tau$. In this work we will suppose that the function $f(w)$
satisfies the following conditions implied by its biological meaning (Figure \ref{f7}):

\begin{equation}\label{m4}
 f(w) > 0  \;\; {\rm for} \;\;  0 \leq w < 1, \;\; f(1) = 0 , \;\; f'(1) > -1 ,
\end{equation}

\begin{equation}\label{m5}
  f(0) > 1 , \; f'(0) > 0,\  f(w) > 1 \;\; {\rm for} \;\; 0 \leq w < w_* ,
\end{equation}
for some $w_* \in (0,1)$. Furthermore,
\begin{equation}\label{m6}
 {\rm equation} \;\;  f(w) = 1-w \;\; {\rm has \; a \; single \; solution} \;\; w_0 \;\; {\rm for} \;\;  0 < w < 1;
 \;\;\; f'(w_0) < -1.
\end{equation}
Hence $f(w) > 1-w$ for $0  \leq w < w_0$ and $f(w) < 1-w$ for $w_0 < w < 1$.
Under these conditions, the function $F(w) = w(1-w-f(w))$ has three zeros: $w=0,$ $ w=w_0,$ $ w=1$. Moreover
$F'(0)<0,$ $ F'(1)<0$.  In the other words,  we consider so-called bistable travelling waves.  We recall here that  travelling wave is a solution of equation (\ref{m0}) having the form $v(x,t) = w(x-ct)$, where the constant
$c$ is the wave speed. Clearly, the wave profile $w(z)$ satisfies  the relation
\begin{equation}\label{m1}
  w'' + cw' + w (1-w - f(w(z+c\tau))) = 0 .
\end{equation}

Equation (\ref{m0}) is  a quasi-linear functional reaction-diffusion equation.  The basic concepts of the general theory of these equations were developed in \cite{MS,wB}.
 In this respect, the delayed reaction-diffusion equations whose reaction term $g$ is either of logistic type (i.e. $g$ is as in (1.1), when  $v$ is not separated multiplicatively from $v_\tau$) or of the Mackey-Glass type (when  $v$ is separated multiplicatively from $v_\tau$, i.e. $g= -kv + b(v_\tau)$) are between the most studied ones, e.g. cf. \cite{GTLMS,TPT}.   It is an interesting point of discussion whether the Mackey-Glass type models reflect more adequately the biological reality than the logistic models, e.g. see \cite[p. 56-58]{Tur} and especially  \cite[Section 1.1]{HD}  for further details. Importantly, in  certain relevant situations both models exhibit similar types of qualitative behaviour of solutions.    It is also worth noting that    the investigation of delayed logistic  models is more  difficult and  technically involved   than the studies of the  Mackey-Glass type systems, precisely because of the multiplicative non-separateness  of $v$ and $v_\tau$.  For example, so far no analytical results on the existence and uniqueness  of bistable waves in  delayed equations which include model  (\ref{m0}) with  {\it non-monotone nonlinear}   response $f$ were available in the literature, cf. \cite{ACR,GSW,ma,FW,sch,SZ,WLR}.

In order to understand what kind of  results can be expected in the delayed case, first we recall
the main existence assertion \cite{VVV} about bistable waves in the classical reaction-diffusion equation without  delay ($\tau=0$):
\begin{equation}\label{m2}
 w'' + cw' + F(w) = 0, \quad  F(0)=F(1)=0, \ F'(0) <0, \ F'(1)<0.
\end{equation}
In this  case, for a unique value of $c$,  there exists a unique (up to translation) monotonically
decreasing solution of problem (\ref{m2}) on the whole axis with the limits
\begin{equation}\label{m3}
  w(-\infty)=1, \;\;\; w(\infty)=0.
\end{equation} 
Furthermore,  in the delayed and monotone case (i.e. when $\tau>0$ and $f(w)$ is decreasing),  existence of solutions for problem (\ref{m1}), (\ref{m3}) was proved in \cite[Theorem 5.8]{WLR}.  In such a case,  equation (\ref{m0}) admits the maximum and comparison
principles, and existence of solutions can be studied using these conventional techniques. See also recent work by Fang and Zhao for a similar result  \cite[Theorem 6.4]{FW} obtained in an abstract setting of monotone bistable semiflows.  However, the properties of equation  (\ref{m1})   change seriously if the function $f(w)$ is not monotonically decreasing.
In this case,
a first attempt  to tackle the aforementioned  existence and uniqueness  problem was recently made  in  \cite{PLOS} .  The idea  of  \cite{PLOS} was to consider a discontinuous piece-wise constant approximation  instead of  the original continuous  function $f$, in this way  simplifying the model  (\ref{m0}) and allowing the use of the phase plane method.

In any event,   the approaches of \cite{PLOS,FW,WLR} are not applicable in the  case of continuous non-monotone $f$, and the question about the existence of travelling waves remains open.

  In this paper, we answer this question affirmatively, by proposing  a different approach based on the construction of the topological degree for an elliptic  operator considered on a subset  containing only bistable waves  monotonically decreasing on the whole real axis. Hence, basically  we are going to use the Leray-Schauder method not
for all solutions of   (\ref{m1}),  (\ref{m3}) but only for monotonically decreasing  ones. It appears that during
a continuous deformation monotone  waves are separated from the non-monotone waves in the sense
that the norm of their difference is uniformly bounded from below by a positive constant, cf. \cite{GTLMS,TPT}. 
This property allows us to construct a domain in the function space which contains all monotone  waves
and which does not contain any non-monotone solution of   (\ref{m1}),  (\ref{m3}). We prove that the value of the topological degree
of the corresponding operator is different from zero in this domain.  Then  the existence of solutions to  (\ref{m1}),  (\ref{m3})  follows.
This procedure was successfully applied in \cite{VVV} to some classes of reaction-diffusion systems while the recent works
\cite{ADV,ADV2,AV,v2011} suggested that it can also be extended to the framework of the theory of  reaction-diffusion equations with delay and bistable nonlinearity.  The key ingredient of this method is the construction of  topological degree for elliptic boundary problems similar to    (\ref{m1}),  (\ref{m3}). Essentially this construction was proposed in \cite{VVV} and recently developed further  in \cite[Chapter 11]{v2011} with \cite{ADV,ADV2,AV}.  It is worth to note that we apply the Leray-Schauder method in a rather direct fashion, without the use of truncation argument, cf. \cite{ACR,BNPR}.  In any case, in the theory of functional reaction-diffusion equation,  the wave profile equations (as (\ref{m1}) or equivalent integral equations) are usually solved either through  the iteration procedure or by means of the Schauder fixed point theorem. These approaches  lead to restrictive  monotonicity assumptions on the delayed term (in some situations,  squeezing technique allows to weaken them and consider  non-monotone delayed terms as well). In this way, the application of the Leray-Shauder method, which is less demanding in regard to the shape and smoothness properties of the delayed nonlinearity, seems to be an interesting new possibility  in this area of research.

Finally, we say a few words about the organisation of the paper.  In the next section,  we analyse briefly the basic properties of  nonlinear operators required for the degree construction: by \cite{v2011}, these operators should be proper and  their Frechet derivatives augmented by the term $-\lambda I$ should be Fredholm operators with  zero index for all $\lambda \geq 0$. Section 3 is devoted to separation
of monotone waves from non-monotone ones while in  Section 4 we establish  a priori estimates of waves. Finally, in Section 5  we prove our main result in this paper, Theorem 5.1: it says that if $C^4-$smooth $f$ satisfies (\ref{m4}), (\ref{m5}),  (\ref{m6})  and $f'(w) <0$ for all $w\in [f^{-1}(1),1]$, then at least one monotone bistable wave  for problem (\ref{m1}),  (\ref{m3}) exists for every fixed delay $\tau \geq 0$.



\setcounter{equation}{0}

\section{Operators and  topological degree}

\medskip

\paragraph{ Fredholm property of  associated linear operators.}

Let $E$ be the H\"older space $C^{2+\alpha}(\mathbb R)$ which consists of functions continuous in $\mathbb R$ together with their second
derivatives and the second derivative satisfies the H\"older condition with the exponent $\alpha \in (0,1)$.
 Recall that the H\"older semi-norm of some function $p: \mathbb R \to \mathbb R$ is defined as
$$
[p]_\alpha:= \sup_{x\not =y}\frac{|p(x)-p(y)|}{|x-y|^\alpha}.
$$
In this section, we will use several times  the following obvious estimate of $[p]_\alpha$ for $p \in C^1(\mathbb R)$:
\begin{equation}\label{Hi}
 [p]_\alpha \leq \max\{2|p|_\infty, |p'|_\infty\}.
 \end{equation}
 Here we use the standard notation  $|p|_\infty= \sup_{x \in\mathbb R }|p(x)|$ for bounded functions $f:\mathbb R  \to \mathbb R$.

\medskip

Similarly, $F$ will denote the space  of H\"older continuous functions with the same exponent $\alpha$.     The norms
in spaces $E$ and $F$ are given by the formulas $\| p \|_E: = |p|_\infty + |p'|_\infty + |p''|_\infty + [p'']_\alpha$ and $\| p \|_F: = |p|_\infty +  [p]_\alpha$, respectively. We also will consider weighted spaces $E_\mu$ and $F_\mu$ defined as follows: $u \in E_\mu$ is and only if $u\mu \in E$,
$\mu(x) =  1+x^2$. Clearly,  $E_\mu$ is a Banach space with the norm $\|f\|_{E_\mu} = \|\mu f\|_E$.
The space $F_\mu$ is defined similarly. The choice of the weight function is not unique.
It can be any positive function with polynomial growth at infinity.    

\medskip

Now, for some fixed real parameter $h$, consider the linear operator $L : E \to F$,
$$ L u = u'' + a(x) u' + b(x) u + d(x) u_h , $$
where $u_h(x) = u(x+h)$. We assume that the coefficients $a(x), b(x)$ and $d(x)$ belong to the space $E$.
Assuming, further, that the coefficients have limits at infinity,
$$ a(x) \to a_\pm , \;\; b(x) \to b_\pm, \;\; d(x) \to d_\pm , \;\;\; x \to \pm \infty , $$
we can introduce the limiting operators
$$ L_\pm u = u'' + a_\pm u' + b_\pm u + d_\pm u_h . $$
Let us recall that the operator $L$ is called normally solvable if its image is closed.
This definition is equivalent to the condition that the equation $Lu=f$ is solvable if and only if
$f$ is orthogonal to all functionals from some closed subspace of the dual space $F^*$.

\medskip

\noindent
{\bf Condition NS.}
The operator $L$ is said to satisfy Condition NS if equations $L_\pm u = 0$ do not have nonzero
bounded solutions.

\medskip

\noindent
{\bf Proposition 2.1.}
{\em The operator $L$ is normally solvable with a finite dimensional kernel if and only if Condition NS
is satisfied.}

\medskip

\noindent
 {\bf Proof.}
The demonstration  of this proposition is similar to the proof of Theorem 2.2 in \cite{ADV} (sufficiency) and
of Theorem 2.3 in \cite{ADV2} (necessity), and therefore it is omitted here.  \hfill $\square$

\medskip

If we substitute $\exp(i\xi)$ in the equations $L_\pm u=0$, then we obtain

\begin{equation}\label{op1}
  -\xi^2 + a_\pm i \xi + b_\pm + d_\pm e^{i \xi h} = 0 .
\end{equation}
Condition NS is satisfied if and only if these equations do not have solutions for any real $\xi$ \cite{VVC}.
Similarly, equations
\begin{equation}\label{op2}
  L_\pm u = \lambda u
\end{equation}
 with $u= \exp(i\xi)$ are equivalent to
\begin{equation}\label{op3}
  \lambda = -\xi^2 + a_\pm i \xi + b_\pm + d_\pm e^{i \xi h} , \;\; \xi \in \mathbb R .
\end{equation}

\noindent
{\bf Theorem 2.2.}
{\em If the curves $\lambda(\xi)$ given by (\ref{op3}) are in the open left-half plane of the
complex plane for all real $\xi$, then the operator   $L- \lambda:E \to F$ with $\lambda \geq 0$ satisfies the Fredholm property, and
its index equals $0$.}

\medskip

\noindent
{\bf Proof.} Note that the operator $L-\lambda$ is normally
solvable with a finite dimensional kernel for all real $\lambda \geq 0$. Therefore its index
is constant for such values of $\lambda$. On the other hand, the operator  $L-\lambda$ is invertible for
$\lambda$ sufficiently large (  to see the latter, it suffices to transform  equation $(L-\lambda)u=f$  into equivalent integral equation and then apply Banach contraction principle  for all large $\lambda >0$).  Therefore its index is $0$. \hfill $\square$

\medskip

Finally, recall  that   the set of complex numbers $\lambda$ for which the operator $L - \lambda $ does not satisfy the Fredholm property is called the essential spectrum of the operator $L$. It is known that a) the essential spectrum of $L$ coincides with the union of two curves given by  (\ref{op3}) and that b) polynomial weight does not change the essential spectrum,  cf. \cite[Chapter 5]{v2011}.    Thus Theorem 2.2 remains valid if we replace in its statement action $L: E \to F$ with $L: E_\mu \to F_\mu$.
\paragraph{Nonlinear operators.} As we have already mentioned in the introductory section,
in this work we study the existence of solutions of the problem
\begin{equation}\label{p1}
  w'' + c w' + w (1-w-f(w(x+c\tau))) = 0 , \;\;\; w(-\infty)=1, \; w(\infty)=0 ,
\end{equation}
where $c$ is an unknown constant which should be chosen to provide the existence of solution.
 The function $f(w)$ is supposed to be $C^4$-smooth and have uniformly bounded derivatives.
In  the case of equation (\ref{p1}), the limiting operators $L_\pm$ introduced in the previous subsection have the following forms:
$$
(L_+u)(x) = u''(x) + cu'(x) +(1-f(0))u(x),  \quad  (L_-u)(x) = u''(x) + cu'(x) - u(x) -f'(1)u(x+c\tau).
$$
Since $f(0) >1$, the operator $L_+$ clearly satisfies the assumption NS. The same property holds
for the operator  $L_-$ (considered with arbitrary $\tau$) because of the inequalities  $-1 < f'(1) \leq 0$. Indeed, it is immediate to see  that the characteristic equation
$z^2 +cz -1 = f'(1)e^{c\tau z}$ associated with $L_-$, can not have pure imaginary solutions.  Furthermore,
the respective curves
$$\lambda_+(\xi) = -\xi^2 + ic \xi + 1-f(0); \quad  \lambda_-(\xi) = -\xi^2 + ic \xi  - 1 -f'(1)e^{ic\tau \xi}  $$
satisfy, for all $\xi \in \mathbb R$, the inequalities
$$\Re \lambda_+(\xi) = -\xi^2  +1-f(0) <0; \quad  \Re \lambda_-(\xi) = -\xi^2  - 1 -f'(1)\cos (c\tau \xi)<0. $$
All the above implies that,  in case of
equation (\ref{p1}),  Condition NS, Proposition 2.1 and Theorem 2.2 can be used without restrictions.

In order to introduce the nonlinear operator corresponding to problem (\ref{p1}), we
set $w(x) = u(x) + \psi(x)$, where $\psi(x)$ is an infinitely differentiable non-increasing function,
$\psi(x) \equiv 1$ for $x \leq 0$, $\psi(x) \equiv 0$ for $x \geq 1$.
Then we consider the  nonlinear operator
$$ A_\tau(u) = (u+\psi)'' +  c(u)(u+\psi)' + (u+\psi) ( 1 - u - \psi - f( u(x+ c(u)\tau) + \psi(x+ c(u)\tau))), $$
$$
c(u)= \ln \sqrt{\int_{ \mathbb R} (u(s)+\psi(s))^2\min\{e^{s},1\}ds}=: 1/2\ln \rho(u),
$$
acting from the space $E_\mu \times \mathbb R$ ($u \in E_\mu, \tau \in \mathbb R$) into the space $F_\mu$.
The main purpose of introducing functional $c(u)$ in the definition of $A_\tau(u)$ instead of considering $c$ as a  real parameter (as in  (\ref{p1}))  is twofold: first, this obliges all solutions of the equation $A_\tau(u)=0,$ $ \ \tau \in [0,\tau_*]$, belong to an open bounded set $D\subset E_\mu$; second, for $\tau=0$, it allows to calculate topological degree $\gamma (A_0,D)$ by removing zero eigenvalues of some associated  linear operator. Since $\gamma (A_0,D)$ was already found in \cite[Chapter 3, \S 3.2 ]{VVV}: $\gamma (A_0,D) =1$, we will use $c(u)$ only to obtain necessary a priori estimates of the wave solutions.

It is easy to see that functional $c:E_\mu \to \mathbb R$ is $C^1$-smooth and
$$
c'(u)h(x) = \frac{1}{\rho(u)}\int_{ \mathbb R} (u(s)+\psi(s))h(s)\min\{e^{s},1\}ds.
$$
Therefore, to prove that the operator $A_\tau(u)$ depends $C^1$-smoothly on $u, \tau$ it suffices to establish
that the nonlinear part of $A_\tau(u)$,  i.e. $\mathfrak{N}(u,\tau):=$
$$
F(u+\psi, u(\cdot+ c(u)\tau) + \psi(\cdot+c(u)\tau)) := (u+\psi) ( 1 - u - \psi - f( u(\cdot+ c(u)\tau) + \psi(\cdot+  c(u)\tau))),
$$
is continuously differentiable.  It is worth to mention that the above formula for $\mathfrak{N}(u,\tau)$ contains state-depending shifts of arguments (i.e. expressions like $u(x+ c(u)\tau)$)  and therefore the differentiability question for $\mathfrak{N}(u,\tau)$  should be handled with certain  care, e.g. see \cite[Section 3]{FH}.  Let us show, for example,  the existence of the Fr\'echet derivative $D_u\mathfrak{N}$. We claim that $[D_u\mathfrak{N}(u,\tau)]h(x) = $
$$
F_1(P_0(x))h(x)+F_2(P_0(x))\left(h(x+c(u)\tau)+  \tau[u'(x+c(u)\tau)+\psi'(x+c(u)\tau)]c'(u)h\right),
$$
where $F_j= F_j(v_1,v_2)$ denotes the partial derivative of $F(v_1,v_2)$ with respect to $v_j$:
$$
F_j(P_0(x)) = F_j(u(x)+\psi(x), u(x+ c(u)\tau) + \psi(x+ c(u)\tau)),
$$
$$
P_s(x)= (u(x)+sh(x)+\psi(x), u(x+c(u+sh)\tau)+ sh(x+c(u+sh)\tau)+ \psi(x+c(u+sh)\tau))\in {\mathbb R}^2.
$$
Indeed, a straightforward  computation shows that
$$
\mathfrak{N}(u+h,\tau)(x)- \mathfrak{N}(u,\tau)(x)- [D_u\mathfrak{N}(u,\tau)]h(x)= R_1(u,h)(x)h(x) + R_2(u,h)(x),
$$
where $R_1(u,h)(x)=\int_0^1(F_1(P_s(x))-F_1(P_0(x)))ds,$
$$
R_2(u,h)(x)= \int_0^1(F_2(P_s(x))-F_2(P_0(x)))\Lambda(x,s)ds+\int_0^1F_2(P_0(x))(\Lambda(x,s)-\Lambda(x,0))ds,
$$
$$
\Lambda(x,s) = \tau(u'(x+c(u+sh)\tau)+ sh'(x+c(u+sh)\tau)+ \psi'(x+c(u+sh)\tau))(c'(u+sh)h)+
$$
$$
h(x+c(u+sh)\tau).
$$
Since $f\in C^4(\mathbb R)$, $u \in E_\mu$,  we find that, for some positive $k_u$ and $\delta_u$
depending only on $u$,
$$
|R_1(u,h)|_\infty \leq k_u|h|_\infty, \
|\mu R_2(u,h)|_\infty \leq k_u\|h\|_{E_\mu}|h|_\infty, \ \ \mbox{for all} \
 \ \|h\|_{E_\mu} \leq \delta_u;
$$
$$
|R_1'(u,h)|_\infty \leq k_u(|h|_\infty+|h'|_\infty), \
|\mu R'_2(u,h)|_\infty \leq k_u\|h\|^{1+\alpha}_{E_\mu}, \ \ \mbox{for all} \
 \ \|h\|_{E_\mu} \leq \delta_u.
$$
In view of (\ref{Hi}), the above estimates imply that $\|R_1(u,h)h+ R_2(u,h)\|_{F_\mu}
 = O(\|h\|^{1+\alpha}_{E_\mu})$. This assures the differentiability  of $\mathfrak{N}(u,\tau)$ with respect to $u$. A similar reasoning also shows that  $D_u\mathfrak{N}$ depends continuously on $u, \tau$  in the operator norm.
 Finally,
$$
[D_uA_\tau(u)]h(x)= h''(x)+c(u)h'(x)+ (c'(u)h)u'(x) + [D_u\mathfrak{N}(u,\tau)]h(x)=
$$
$$
h''(x)+a(x)h'(x)+ b(x)h(x) + d(x)h(x+ {c(u)\tau})+ (Kh)(x) = (Lh)(x) + (Kh)(x),
$$
where $a(x) = c(u), \ b(x)= F_1(P_0(x)), \ d(x)= F_2(P_0(x))$ and $K: E_\mu \to F_\mu$ defined by
$$
(Kh)(x) = (c'(u)h)\left[u'(x)+\tau u'(x+c(u)\tau)+ \tau \psi'(x+c(u)\tau)\right]
$$
is one-dimensional linear  operator. Since finite dimensional perturbations of the Fredholm operator
does not change its index,  we obtain the following version of Theorem 2.2:
\medskip

\noindent
{\bf Theorem 2.3.}
{\em If the curves $\lambda(\xi)$ given by (\ref{op3}) are in the open left-half plane of the
complex plane for all real $\xi$, then the operator   $D_uA_\tau- \lambda: E_\mu \to F_\mu$ with $\lambda \geq 0$ satisfies the Fredholm property, and
its index equals $0$.  Moreover, there exists $C_0, \lambda_0>0$ such that
\begin{equation}\label{inv}
\|(D_uA_\tau- \lambda)^{-1}\| \leq C_0 \ \mbox{for all} \ \lambda > \lambda_0, \ \tau \in [0, \tau_*].
\end{equation}

\noindent {\bf Proof.} In view of the above said, we only have to prove inequality (\ref{inv}).  We have that
$$(D_uA_\tau)h(x)- \lambda^2 h(x) = h''(x)- \lambda^2 h(x)+a(x)h'(x)+ b(x)h(x) + d(x)h(x+ c\tau)+ \phi(x)l(h),$$
where linear functional $l(h):= c'(u)h,\ l:E_\mu \to F_\mu,$ is continuous and $\phi, \phi' \in F_\mu$.  Set $v= \mu h$, then
$$
h''(x)- \lambda^2 h(x)+a(x)h'(x)+ b(x)h(x) + d(x)h(x+ c\tau)+ \phi(x)l(h) = g(x)
$$
if and only if
\begin{equation}\label{fv}
v''(x)- \lambda^2 v(x)+a_1(x)v'(x)+ b_1(x)v(x) + d_1(x)v(x+ c\tau)+ \phi_1(x)l_1(v) = g_1(x),
\end{equation}
where $g_1(x)=\mu(x)g(x),\ \phi_1(x)=\mu(x)\phi(x),  \ l_1(v) = l(v/\mu)$ and
$$
a_1(x)= a(x)-\frac{2\mu'(x)}{\mu(x)}, \ b_1(x)= b(x) +\frac{2(\mu'(x))^2}{(\mu(x))^2}- \frac{\mu''(x)}{\mu(x)} - a(x)\frac{\mu'(x)}{\mu(x)}, d_1(x) = d(x)\frac{\mu(x)}{\mu(x+c\tau)}.
$$

Each bounded solution of equation (\ref{fv}) should satisfy the integral equation
\begin{equation}\label{int1}
v = \frac{1}{2\lambda}(M_+Tv - M_+g_1),
\end{equation}
where
$$
M_\pm w(t) = \left(\pm \int_{-\infty}^te^{-\lambda(t-s)}w(s)ds+ \int^{+\infty}_te^{\lambda(t-s)}w(s)ds\right),
$$
$$
T(x, v(x),v'(x))= a_1(x)v'(x)+ b_1(x)v(x) + d_1(x)v(x+ c\tau)+ \phi_1(x)l_1(v).
$$
After differentiating (\ref{int1}), we find that
\begin{equation}\label{int2}
v' = \frac{1}{2}(M_-Tv - M_-g_1).
\end{equation}
Consider $M: C({\mathbb R}, {\mathbb R}^2) \to C({\mathbb R}, {\mathbb R}^2)$ defined by
$$M(v,w) = \frac 1 2 (\frac{1}{\lambda}M_+Tv, M_-Tv),$$ it is immediate to see that $\|M|| \leq K\lambda^{-1}$, where $K = 1+ |a_1|_\infty+|b_1|_\infty+|d_1|_\infty + |\phi_1|_\infty \|l_1\|$ does not depend on $\tau$.  Therefore system (\ref{int1}), (\ref{int2})  has a unique solution $(v, v') \in  C({\mathbb R}, {\mathbb R}^2)$ once $\lambda > K =\lambda_*$. Moreover,
$$
|(v,v')|_\infty \leq 2 |(M_+g_1,M_-g_1)|_\infty \leq \frac{4}{\lambda}|g_1|_\infty \ \mbox{for all} \ \lambda > 2\lambda_*.
$$
Using this inequality in equation (\ref{int1}) again, we can improve the estimate for $|v|_\infty$ as follows:
$
|v|_\infty \leq K_1 \lambda^{-2}|g_1|_\infty,  \ [v]_\alpha \leq \lambda^{-2}([g_1]_\alpha+ [Tv]_\alpha), \  \lambda > 2\lambda_*.
$ Then (\ref{fv}) yields that $|v''|_\infty \leq K_2|g_1|_\infty$ for all $\lambda > 2\lambda_*$. Therefore, using  (\ref{Hi}), we also find that $ [Tv]_\alpha \leq K_3 [g_1]_\infty$. Hence,  $ [v]_\alpha \leq K_4\lambda^{-2}\|g\|_{F_\mu}$ for all $\lambda > 2\lambda_*$.
Here $K_j, j =1,2,3,4,$ stand for some universal constants.
In this way,
$$
[v'']_\alpha \leq \lambda^2 [v]_\alpha+[a_1v']_\alpha+ [b_1v]_\alpha + [d_1v]_\alpha +  [\phi_1]_\alpha |l_1(v)| + [g_1]_\alpha \leq K_5\|g\|_{F_\mu},  \ \lambda > 2\lambda_*,
$$
so that
$$\|(D_uA_\tau- \lambda)^{-1}g\|_{E_\mu} =  \|h\|_{E_\mu} =   \| v\|_{E} \leq  K_5\|g\|_{F_\mu},   \ \lambda > 2\lambda_*.$$
The latter  estimate  completes the proof of Theorem 2.3.
 \hfill $\square$

\medskip

The nonlinear operator $A_\tau$ has another useful  property: it is a proper operator} in the sense that the inverse image of compact
sets is compact in any bounded closed set:
\medskip

\noindent
{\bf Theorem 2.4.}
{\em Let the condition NS
be satisfied for $D_{u_0}A_{\tau_0}$.  Assume $A_{\tau_n}(u_n) =g_n$ for  converging sequences
of elements $g_n\in F_\mu,$  $\tau_n \geq 0$. If, in addition,  sequence $\{u_n\}$ is bounded in $E_\mu$, then it has a  subsequence converging in $E_\mu$. }

\medskip

\noindent
{\bf Proof.} Suppose that $|u_n|_{E_\mu} \leq K$. Since the inclusion $E_\mu \subset C^{2}(\mathbb R)$ is compact, there exists a subsequence $u_{n_j}$ converging in $C^{2}(\mathbb R)$ to some element $u_0$.
Clearly, $v_0= \mu u_0 \in C^{2,\alpha}$ and  $|v_0|_{E} \leq K$. Without loss of generality, we can assume that $u_n \to u_0$, then  $c_n:= c(u_n) \to c_0 = c(u_0)$ and also $v_n: = \mu u_n \to v_0$ uniformly on compact subsets of $\mathbb R$. The same type of convergence holds for the first and second derivates of $v_n$.
To prove that $v_n: = \mu u_n \to v_0$ uniformly on $\mathbb R$, we consider the
following relation:
\begin{equation} \label{tak}
\mu\left[A_{\tau_n}\left(\frac{v_n}{\mu}\right)- A_{\tau_0}\left(\frac{v_0}{\mu}\right)\right] = \mu(g_n-g_0).
\end{equation}
Suppose for a moment that $V_n= v_n-v_0$ does not converge,  uniformly  on $\mathbb R$, to the zero function.
Then there exist positive $\epsilon_0$ and  some sequence $x_n$ such that $|V_n(x_n)| \geq \epsilon_0$ for all $n$.  Since sequence $\{x_n\}$  can not be bounded, we may suppose first that  $x_n \to + \infty$.  Now, since $W_n(x)= V_n(x+x_n)$ satisfies $|W_n|_E \leq 2K$,  without restricting generality, we can assume that
$W_n(x)$, together with their first and second derivatives, converges uniformly on compact subsets of $\mathbb R$. Let $W_0(x)$ be the limit function for  $\{W_n(x)\}$, then $|W_0(x)| \leq 2K,$ $|W_0'(x)| \leq 2K,$ $x \in \mathbb R,$ and $|W_0(0)| \geq \epsilon_0$.  After taking limit $n \to +\infty$ in (\ref{tak}), since $\psi(x+x_n) \to 0$ uniformly on compact sets,  it is easy to find that
$$
W_0''(x)+c_0W_0'(x) + (1-f(0))W_0(x) =0,
$$
where $(1-f(0))W_0(x) = \lim_{n \to +\infty} R_n(x+x_n)$ and
$$
R_n(x) = (v_n(x)+\mu(x)\psi(x)) ( 1 - u_n(x) - \psi(x) - f( u_n(x+c_n\tau_n) + \psi(x+c_n\tau_n)))-
$$
$$
(v_0(x)+\mu(x)\psi(x)) ( 1 - u_0(x) - \psi(x) - f( u_0(x+c_0\tau_0) + \psi(x+c_0\tau_0))).
$$
Since the only bounded solution of the latter differential equation is $W_0(x) \equiv 0$, we arrive to a contradiction
(recall that $W_0(0)\not=0$).

Now, if $x_n\to -\infty,$ then $\psi(x+x_n) \to 1$ uniformly on compact sets and therefore
$$
 \lim_{n \to +\infty} R_n(x+x_n) = -W_0(x)(1+f(1))-
 $$
 $$
 \lim_{n \to +\infty} \tilde W_n(x) \left (\frac{f( u_n(x+x_n+c_n\tau_n) +1))- f( u_0(x+x_n+c_0\tau_0) +1))}{u_n(x+x_n+c_n\tau_n)- u_0(x+x_n+c_0\tau_0)}\right) =
$$
$$
 -W_0(x)(1+f(1))-  W_0(x+c_0\tau_0)f'(1),
$$
where $\lim \tilde W_n(x) = W_0(x+ c_0\tau_0)$ because of
$$
 \tilde W_n(x):= \mu(x+x_n)(u_n(x+x_n+c_n\tau_n)- u_0(x+x_n+c_0\tau_0)) =
$$
$$
W_n(x+c_n\tau_n)+  \mu(x+x_n)(u_0(x+x_n+c_n\tau_n)- u_0(x+x_n+c_0\tau_0)) =$$
$$ W_n(x+c_n\tau_n)+ \mu(x+x_n)u_0'(x+x_n+ \theta_n) (c_n\tau_n - c_0\tau_0).
$$
Observe here that $u_n(x+x_n+c_n\tau_n), u_0(x+x_n+c_0\tau_0)  \to 0 \ \mbox{ as } \ n\to +\infty.$

Hence, $W_0(x)$ satisfies the functional differential equation
$$
W_0''(x)+c_0W_0'(x) -  W_0(x)(1+f(1))-  W_0(x+c_0\tau_0)f'(1) =0,
$$
which in view of condition NS can have only zero bounded solution.  The obtained  contradiction shows  that  $W_n(x) $ converges to $0$  uniformly on $\mathbb R$ (so that
$v_n=\mu u_n \to v_0=\mu u_0$ uniformly on $\mathbb R$).

Similarly, if $V'_n(x)$ does not converge uniformly on $\mathbb R$ to $0$,  then $|V'_n(x_n)| \geq \epsilon_0>0$ for some $\epsilon_0$ and $\{x_n\}$.  Considering $W_n(x)=V_n(x+x_n)$,  we obtain from the the previous part of the proof  that  $W_n(x) \to 0, \ n \to +\infty$ (uniformly on $\mathbb R$). Therefore, since $W'_n(x)$ converges, uniformly on compact sets,  to some continuous function $W_*(x)$, we conclude that $W_*(x) \equiv 0$. Clearly, this contradicts
to the inequalities $|W'_n(0)|\geq \epsilon_0, \ n =1,2,\dots$ Hence, we can conclude that  $V'_n(x)$  converges to $0$ uniformly on $\mathbb R$. The proof of the uniform (on $\mathbb R$) convergence $V''_n(x) \to 0$  is based on  the same argument.  Finally, in order to estimate the H\"older semi-norm
$
[V_n'']_\alpha$,
it suffices to use the following equivalent form of (\ref{tak}):
\small{
$$
V_n''(x)= \mu(x)(g_n(x)-g_0(x))+ V_n'(x)\left(\frac{2\mu'(x)}{\mu(x)}-c_0\right) + (c_0-c_n)\left(\psi'(x)\mu(x) -v_n\frac{\mu'(x)}{\mu(x)}+ v_n'(x)\right)+
$$
$$
V_n(x)\left(\frac{\mu''(x)}{\mu(x)} -\frac{2(\mu'(x))^2}{\mu^2(x)}+2\psi(x)-1+c_0\frac{\mu'(x)}{\mu(x)}-u_0(x)-u_n(x)+f(u_0(x+c_0\tau_0)+\psi(x+c_0\tau_0))\right)
$$
 $$+(\mu(x)+v_n(x))\left(f(u_n(x+c_n\tau_n)+\psi(x+c_n\tau_n))-f(u_0(x+c_0\tau_0)+\psi(x+c_0\tau_0))\right).
$$}
It is easy to see from this representation (directly leading to a Schauder type interior estimate) that $\lim_{n\to +\infty} [V_n'']_\alpha =0$. For instance, we can estimate the H\"older semi-norm of the third line in the above expression by using inequality (\ref{Hi})  and the uniform on $\mathbb R$ convergences
$$
\mu(x)\left((u_n(x+c_n\tau_n)+\psi(x+c_n\tau_n))-(u_0(x+c_0\tau_0)+\psi(x+c_0\tau_0))\right) \to 0, \ n \to +\infty,
$$
$$
\mu(x)\left((u_n'(x+c_n\tau_n)+\psi'(x+c_n\tau_n))-(u_0'(x+c_0\tau_0)+\psi'(x+c_0\tau_0))\right) \to 0,\ n \to +\infty.
$$
This completes the proof of Theorem 2.4.
 \hfill $\square$

\medskip

Theorems 2.3 and 2.4 show that, according to the  abstract theory developed in \cite[\S 3.3 of Chapter 11]{v2011}, the Leray-Schauder type topological degree can defined for the operator $A_\tau$. Indeed, in terms of \cite{v2011},  Theorems 2.3 and 2.4 imply that, for a fixed $\tau\geq 0$, the operator  $A_\tau$ belongs to the class ${\mathbf F}$ while the family of operators  $A_\tau$  with $\tau \in [0,\tau_0]$ belongs to the class ${\mathbf H}$.

\setcounter{equation}{0}

\section{Monotonicity of solutions}

In this section, we will consider  the wave profile equation
\begin{equation}\label{b1}
  w'' + cw' + w (1-w - f(w(x+c\tau)))=0
\end{equation}
with the boundary conditions
\begin{equation}\label{b4}
  w(-\infty) = 1 , \;\;\; w(+\infty) = 0 .
\end{equation}

\medskip

\noindent
{\bf Lemma 3.1.}
{\em Suppose that $c \geq 0$, $f'(w) < 0$ for $w_0 \leq w < 1$. If solution $w(x)$ of problem
(\ref{b1}), (\ref{b4}) satisfies condition $w'(x) \leq 0$ for all $x \in \mathbb R$, then
$w'(x) < 0$, $x \in \mathbb R$.}

 \medskip

\noindent
{\bf Proof.}
Suppose that the assertion of the lemma does not hold and $w'(x_0)=0$ for some $x_0$.
Then $w''(x_0)=0$ and from the equation (\ref{b1}) we obtain the equality
$$ w(x_0) (1 - w(x_0) - f(w(x_0+c\tau))) = 0 . $$
If $w(x_0)=0$, then by virtue of the uniqueness of solution $w(x) \equiv 0$, and we obtain a contradiction
with (\ref{b4}). Hence
\begin{equation}\label{b5}
  1 - w(x_0) = f(w(x_0+c\tau)) .
\end{equation}
Since $c \tau >0$, then $w(x_0+c\tau) \leq w(x_0)$.
Set $w_2 = w(x_0+c\tau)$. Then $w_2 \leq w(x_0)$. From (\ref{b5}) we get

\begin{equation}\label{b6}
  1 - w_2 \geq f(w_2) .
\end{equation}
Suppose that $w(x_0) < w_0$.
Since $1 < f(0)$, that is $1-w < f(w)$ for $w=0$, then by virtue of (\ref{b6}), equation $f(w)=1-w$ has a solution in the interval
$0<w\leq w_2 \; (\leq w(x_0) <w_0)$. This conclusion contradicts the assumption on the function $f(w)$.
Hence $w(x_0) \geq w_0$.

Next, we show that $w_2 \geq w_0$. Indeed, suppose that $w_2 < w_0$. Then (\ref{b6})
contradicts the assumption that $1-w < f(w)$ for $0 \leq w < w_0$.
Thus, $f'(w_2) < 0$.

Set $u(x) = - w'(x)$. Differentiating equation (\ref{b1}), we obtain
\begin{equation}\label{b7}
  u'' + c u' + a(x)u + b(x) = 0 ,
\end{equation}
where
$$ a(x) = k (1-2 w(x)-f(w(x+c\tau))) , \;\;\; b(x) = -  w(x) f'(w(x+c\tau)) u(x+c\tau) . $$
Let us recall that $u(x) \geq 0$ for all $x$,
$u(x_0)=0$, $u'(x_0)=0$, $f'(w(x_0+c\tau))<0$. Since the function $w(x)$ satisfies (\ref{b4}), then
$u(x) \not\equiv 0$.

Let $I_0$ be the maximal interval containing the point $x=x_0$ and such that $u(x) = 0$ for all $x \in I_0$.
Similar to the arguments presented above we can verify that $f'(w(x+c\tau))<0$ for $x \in I_0$.
We take an interval $I$ slightly larger than $I_0$ such that  $f'(w(x+c\tau))<0$ and $u(x) \not\equiv 0$ for $x \in I$.
If this set is reduced to a single point $x=x_0$, then it is a sufficiently small
interval around this point. Since $b(x) \geq 0$ in this interval, $u(x) \geq 0$ and not identically $0$,
then we obtain a contradiction with the maximum principle for equation (\ref{b7}). \hfill  $\Box$

\medskip

\noindent
{\bf Lemma 3.2.}
{\em Suppose that $c \leq 0$, $f'(w) < 0$ for $w_* \leq w \leq w_0$, where $f(w_*)=1$ and
$f(w)>1$ for $0 \leq w < w_*$. If solution $w(x)$ of problem
(\ref{b1}), (\ref{b4}) satisfies condition $w'(x) \leq 0$ for all $x \in \mathbb R$, then
$w'(x) < 0$, $x \in \mathbb R$.}

 \medskip

\noindent
{\bf Proof.}
Suppose that the assertion of the lemma does not hold and $w'(x_0)=0$ for some $x_0$.
Then we obtain equality (\ref{b5}). Since $c\tau \leq 0$, then  $w(x_0+c\tau) \geq w(x_0)$.
Set $w_2 = w(x_0+c\tau)$. It follows from (\ref{b5}) that
$$ 1-w_2 \leq f(w_2) . $$
Hence $w_2 \leq w_0$.
From (\ref{b5}) we get that $f(w(x_0+c\tau)) \leq 1$. Therefore, $w(x_0+c\tau) \geq w_*$.
Thus, $f'(w(x_0+c\tau)) < 0$ and we proceed with equation (\ref{b7}) as in the proof of Lemma 3.1.
\hfill  \hfill $\Box$

\bigskip

\noindent
{\bf Lemma 3.3.}
 Let $w_n(x)$ be solutions of problem (\ref{b1}), (\ref{b4}) for some $\tau_n \in [0, \tau_*]$, $c=c_n$, $|c_n| \leq c_*$ (for some $c_*$),
$n=1,2,...$.
Suppose that $w_0'(x) < 0$ for all $x \in \mathbb R$ and $w_n(x) \to w_0(x)$ in $C^1(\mathbb R)$,
$c_n \to c_0, \tau_n \to \tau_0$.
If $f'(1) > -1$, then there exists $x=x_0$ such that $w_n'(x) < 0$ for $x \leq x_0$ and
$n$ sufficiently large.

\medskip

\noindent
{\bf Proof.}
Let $\epsilon>0$ be such that
\begin{equation}\label{z2}
  1-w > f(w) , \; 1-\epsilon \leq w < 1 ; \;\;\;
  1-w < f(w) , \; 1 < w < 1 + \epsilon .
\end{equation}
We choose such $x_0$ that $w_0(x) > 1 - \epsilon/2$ for all $x \leq x_0$.
Then for $n$ sufficiently large
$$ 1 - \epsilon < w_n(x) < 1 + \epsilon , \;\; x \leq x_0 . $$
Denote by $M$ a positive constant such that $|c_n \tau_n | < M$ for all $n$.
Then for all $n$ sufficiently large, $w_n(x) < 1$ and $w_n'(x) < 0$ for $x_0 - M \leq x \leq x_0$.

Suppose that the assertion of the lemma does not hold. 
We consider two cases: $w_n(x) > 1$ for some $x \leq x_0$ and $w_n(x) \leq 1$ for all $x \leq x_0$.
In the first case, since $w_n(x) \to 1$ as $x \to -\infty$, then there is a global maximum of this function for $x \leq x_0$:

$$ w_n'(x_n) = 0 , \;\; w_n(x) \leq w_n(x_n) , \; x \leq x_0 , \;\; w_n(x_n) > 1 . $$
We have:
$$ 1 - w_n(x_n) < f(w_n(x_n)) \leq f(w_n(x_n + c_n \tau_n)) . $$
Hence $w_n''(x_n) > 0$ and we obtain a contradiction.

Suppose now that $w_n(x) \leq 1$ for $x \leq x_0$ and all $n$ sufficiently large.
If $w_n'(x) \leq 0$ for $x \leq x_0$ and $w_n'(x_n)=0$ for some $x_n < x_0$, then we obtain a contradiction
with Lemma 3.1. Therefore, if the assertion of this lemma is not satisfied, then the function $w_n(x)$
has a minimum for $x < x_0$. Suppose that there exists the most right minimum $x_*$ of this function for $x \leq x_0$.
By virtue of the construction above, $x_* < x_0 - M$.
Since $w_n''(x_*) \geq 0$, then we conclude from the equation that
$$ 1 - w_n(x_*) \leq f(w_n(x_* + c_n \tau_n)) . $$
It follows from condition (\ref{z2}) that
$$ 1 - w_n(x_*+ c_n \tau_n) > f(w_n(x_* + c_n \tau_n)) . $$
Therefore
\begin{equation}\label{z3}
  w_n(x_* + c_n \tau_n) \leq w_n(x_*) .
\end{equation}
If $c_n < 0$, then there is another minimum $x_{**} < x_*$ of this function, and
$w_n(x_{**}) \leq w_n(x_*)$. Repeating the same arguments, we will obtain a sequence of minima
going to $-\infty$. This contradicts the convergence $w_n(x) \to 1$ as $x \to -\infty$.

If $c_n > 0$, then there is a single maximum $x^*$ of this function in the interval $x_* < x^* < x_0$ since
$x_*$ is the most right minimum and $w_n'(x_0) < 0$.
Suppose first that $x_* + c_n\tau_n \leq x^*$. Then  $w_n(x_* + c_n \tau_n) > w_n(x_*)$, and we obtain a
contradiction with (\ref{z3}). Let now $x_* + c_n\tau_n > x^*$. Since $w_n(x)$ is decreasing for $x > x^*$, then
$$  w_n(x^* + c_n \tau_n) <  w_n(x_* + c_n \tau_n)  $$
and
$$ 1- w_n(x^*) <  1 - w_n(x_*) \leq f(w_n(x_* + c_n \tau_n)) < f( w_n(x^* + c_n \tau_n)) . $$
Therefore, $w_n''(x^*) > 0$ and we obtain a contradiction since $x^*$ is a point of maximum.

If the most right minimum does not exist and there is a sequence of extrema converging to some $\hat x$,
then it is sufficient to take a minimum sufficiently close to $\hat x$ and to repeat similar arguments as above.
Let us also note that for $c_n<0$ it is not necessary to take the most right minimum.
Finally, if $c_n=0$, then we obtain the equation without delay for which the assertion of the lemma is known \cite{VVV}. \hfill  $\Box$

\medskip

\noindent
{\bf Lemma 3.4.}
{\em Let $w_n(x)$ be solutions of problem (\ref{b1}), (\ref{b4}) for some $\tau_n \in [0, \tau_*]$, $c=c_n$, $|c_n| \leq c_*$ (for some $c_*$),
$n=1,2,...$. Suppose that $\tau_n \to \tau_0$, $c_n \to c_0$ and $w_n(x) \to w_0(x)$ in $C^1(\mathbb R)$, where
$w_0(x)$ is a solution of problem (\ref{b1}), (\ref{b4}) for $\tau=\tau_0$, $c=c_0$.
If $w_0'(x) < 0$ for all $x \in \mathbb R$, then $w_n'(x) < 0$ for all $x \in \mathbb R$
and $n$ sufficiently large.}

\medskip

\noindent
{\bf Proof.}
Suppose that the assertion of the lemma does not hold, and there is a sequence $x_n$ such that
$w_n'(x_n)=0$. If this sequence is bounded, then we can choose a convergent subsequence, $x_{n_k} \to x_0$.
Then $w_0'(x_0)=0$, and we obtain a contradiction with the assumption of the lemma.

Consider next the case where $x_n \to \infty$. Let $x=x_*$ be the solution of the equation $w_0(x) = w_* - \epsilon$
for some $\epsilon>0$ sufficiently small. Let us recall that $f(w) > 1$ for $0 \leq w < w_*$.
Then for all $n$ sufficiently large and for all $x \geq x_*$, $w_n(x) \leq w_*$. Hence
$$ f(w_n(x + c\tau)) > 1 \;\; {\rm for} \;\; x \geq x^*= x_* + c_* \tau_* $$
and for $n$ sufficiently large since
$$ x+c\tau \geq x^* + c\tau = (x_* + c_*\tau_*) + c\tau \geq x_* . $$
We use here the assumption that $|c\tau| \leq c_*\tau_*$.
Furthermore, $w_n(x^*) \to w_0(x^*)>0$, $w_n'(x^*) \to w_0'(x^*) < 0$.

Let $x_n > x^*$. If $w_n(x_n) > 0$, then $w_n''(x_n) = - w_n(x_n) (1-w_n(x_n) - f(w_n(x_n+c\tau))) > 0$.
Hence any positive extremum is a minimum, and the function $w_n(x)$ cannot converge to $0$ at infinity.

Suppose now that $w_n(x_n) < 0$. Since $w_n(x) \to 0$ as $x \to \infty$, without loss of generality we
can assume that $x_n$ is a global minimum of the function $w_n(x)$. 
Hence $w_n(x_n) \leq w_n(x_n + c\tau)$. Therefore,
\begin{equation}\label{z1}
  1 - w_n(x_n) < f(w_n(x_n)) \leq f(w_n(x_n + c\tau)) .
\end{equation}
The first inequality in (\ref{z1}) holds since $f(w) > 1-w$ in some neighborhood of $w=0$
(including small negative $w$). The second inequality in (\ref{z1}) takes place
for sufficiently small in absolute value $w_n(x_n)$ and $w_n(x_n + c\tau)$ since
$f'(0)>0$.
Thus, $w_n''(x_n) = -  w_n(x_n) (1-w_n(x_n) - f(w_n(x_n+c\tau))) < 0$.
Therefore $x_n$ is a point of maximum, and it cannot be the global minimum, as supposed.

It remains to note that convergence $x_n \to -\infty$ cannot hold due to Lemma 3.3.
\hfill $\Box$

\setcounter{equation}{0}

\section{A priori estimates}

\subsection{Estimate of the wave speed}

Let us introduce functions $f_0(w)$ and $f_1(w)$ such that
\begin{equation}\label{s1}
  f_0(w) \leq f(w) \leq f_1(w) , \;\; f_0'(w) \leq 0, \;\; f_1'(w) \leq 0, \;\;\; 0 \leq w \leq 1,
\end{equation}
$f_0(0) > 1, f_1(0) > 1$, and equations
$$ 1-w=f_0(w) , \;\;\; 1-w = f_1(w) $$
have unique solutions in the interval $0 \leq w < 1$. Then the problem
\begin{equation}\label{s2}
  w'' + c w' + w (1-w-f_0(w)) = 0 , \;\;\; w(-\infty)=1, \;\; w(\infty)=0
\end{equation}
(without delay)
has a unique solution (up to translation in space) $w_0(x)$ for a unique value $c=c_0$.
Similarly, the problem
\begin{equation}\label{s3}
  w'' + c w' + w (1-w-f_1(w)) = 0 , \;\;\; w(-\infty)=1, \;\; w(\infty)=0
\end{equation}
(without delay)
has a unique solution (up to translation in space) $w_1(x)$ for a unique value $c=c_1$.

\medskip

\noindent
{\bf Lemma 4.1.}
{\em If there exists a monotonically decreasing solution $w(x)$ of problem (\ref{b1}), (\ref{b4}) for some $c>0$, then
$c \leq c_0$.}

\medskip

\noindent
{\bf Proof.}
Since $c\tau > 0$ and $w(x)$ is a decreasing function, then $w(x+c\tau) < w(x)$. Hence

\begin{equation}\label{s4}
  f(w(x+c\tau)) \geq f_0(w(x+c\tau)) > f_0(w(x)) .
\end{equation}
Consider the Cauchy problem

\begin{equation}\label{s5}
  \frac{\partial u}{\partial t} = \frac{\partial^2 u}{\partial x^2} + c
  \frac{\partial u}{\partial x} + u (1-u-f_0(u))
\end{equation}
(without delay) with the initial condition
\begin{equation}\label{s6}
  u(x,0) = w(x) .
\end{equation}
Taking into account equation (\ref{b1}) and inequality (\ref{s4}), we obtain:
$$ w'' + cw' + w(1-w-f_0(w)) = w'' + c w' + w (1-w - f(w(x+c\tau))) + $$ $$ w (f(w(x+c\tau)) - f_0(w(x))) > 0 ,
\;\;\; x \in \mathbb R . $$
Hence $w(x)$ is a lower function, and solution $u(x,t)$ of problem (\ref{s5}), (\ref{s6}) is monotonically increasing with respect to $t$
for each $x$.

On the other hand, by virtue of global stability of monotone waves for the bistable equation,
$u(x,t) \to w_0(x + (c-c_0)t)$ as $t \to \infty$ uniformly on the whole axis. From this convergence we can conclude
that $c \leq c_0$. Indeed, if $c > c_0$, then for each $x$ fixed $w_0(x + (c-c_0)t) \to 0$ as $t \to \infty$.
However, $u(x,t) \geq w(x)$ for all $x$ and $t$. This contradiction completes the proof of the lemma.
\hfill  $\Box$

\medskip

\noindent
{\bf Lemma 4.2.}
{\em If there exists a monotonically decreasing solution $w(x)$ of problem (\ref{b1}), (\ref{b4}) for some $c<0$, then
$c \geq c_1$.}

\medskip

\noindent
{\bf Proof.}
Since $c\tau < 0$ and $w(x)$ is a decreasing function, then $w(x+c\tau) > w(x)$. Hence

\begin{equation}\label{s7}
  f(w(x+c\tau)) \leq f_1(w(x+c\tau)) < f_1(w(x)) .
\end{equation}
Consider the Cauchy problem

\begin{equation}\label{s8}
  \frac{\partial u}{\partial t} = \frac{\partial^2 u}{\partial x^2} + c
  \frac{\partial u}{\partial x} + u (1-u-f_1(u))
\end{equation}
(without delay) with the initial condition
\begin{equation}\label{s9}
  u(x,0) = w(x) .
\end{equation}
Taking into account equation (\ref{b1}) and inequality (\ref{s7}), we obtain:
$$ w'' + cw' + w(1-w-f_1(w)) = w'' + c w' + w (1-w - f(w(x+c\tau))) + $$ $$ w (f(w(x+c\tau)) - f_1(w(x))) < 0 ,
\;\;\; x \in \mathbb R . $$
Hence $w(x)$ is an upper function, and solution $u(x,t)$ of problem (\ref{s5}), (\ref{s6}) is monotonically decreasing with respect to $t$
for each $x$.

On the other hand, by virtue of global stability of monotone waves for the bistable non-delayed equation,
$u(x,t) \to w_1(x + (c-c_1)t)$ as $t \to \infty$ uniformly on the whole axis. From this convergence we can conclude
that $c \geq c_1$. Indeed, if $c < c_1$, then for each $x$ fixed $w_0(x + (c-c_1)t) \to 1$ as $t \to \infty$.
However, $u(x,t) \leq w(x)$ for all $x$ and $t$. This contradiction completes the proof of the lemma.
\hfill  $\Box$

\medskip

From the last two lemmas we obtain the following estimate for the wave speed:
\begin{equation}\label{s10}
 |c| \leq c_*: = \max\{|c_0|, |c_1|\}.
\end{equation}
We note that $c \geq 0$ if $c_1 \geq 0$ and $c \leq 0$ if $c_0 \leq 0$.


\subsection{Estimates of solutions}

In this section, we will use repeatedly  the following simple observation:

\medskip

\noindent
{\bf Lemma 4.3.}
{\em Let $w(x)$ be a monotonically decreasing solution of equation (\ref{b1}) with the limits
$w(-\infty)=w_0, w(+\infty)=0$ at infinity. Then $c<0$. If $w(-\infty)=1, w(+\infty)=w_0$, then $c>0$.}

\medskip

\noindent
{\bf Proof.}  Consider a decreasing wave $w(x)$ connecting $w_0$ and $0$ and suppose that
$c \geq 0$. Then $0 \leq w(x+c\tau) \leq w(x) \leq w_0$ and therefore $1-w(x) \leq f(w(x+c\tau))$ so that
$$
w''(x) = -cw'(x)  - w(x) (1-w(x) - f(w(x+c\tau)))  \geq 0.
$$
Since a convex function cannot connect two final equilibria, we have got a contradiction. Thus $c <0$.

Similarly, suppose that some decreasing wave $w(x)$ connects $1$ and $w_0$ with the speed
$c \leq 0$. Then $1 \geq w(x+c\tau) \geq w(x) \geq w_0$ and therefore $1-w(x) \geq f(w(x+c\tau))$, see Fig. 1.
In this way, we will get a contradiction again:
$$
w''(x) = -cw'(x)  - w(x) (1-w(x) - f(w(x+c\tau)))  \leq 0.
$$
This completes the proof of Lemma 4.3.
 \hfill $\Box$

\medskip

We can now estimate the weighted norm of wave profiles:

\medskip

\noindent
{\bf Lemma 4.4.}
{\em Let $w_\tau(x)$ be a monotonically decreasing solutions of problem (\ref{b1}), (\ref{b4})
for possibly different values of $\tau \in [0,\tau_*]$  such that $w_\tau(0)\in [\alpha, \beta] \in (0,1)$ for some fixed $\alpha, \beta$ and all admissible $\tau \in [0,\tau_*]$. Then   there exists a positive constant $M$ independent of
$\tau$ such that
\begin{equation}\label{s11}
  \sup_{x \in \mathbb R} |w'_\tau(x) - \psi'(x)| \mu(x) +   \sup_{x \in \mathbb R} |w_\tau(x) - \psi(x)| \mu(x) \leq M.
\end{equation}
}
\noindent
{\bf Proof.} Let $c_*$ be as defined in (\ref{s10}).  Denote by $x_1(\tau)$ the solution of the equation $w_\tau(x)=w_1$ and by $x_2(\tau)$ the solution of the
equation $w_\tau(x)=w_2$ where  $0< w_2 < \min\{w_0, \alpha\} \leq \max\{w_0, \beta\}< w_1 <1$ are some fixed numbers.

Clearly,  $x_1(\tau)< 0, \ x_2(\tau) >0$. Moreover,
we claim  that the difference $x_1(\tau)-x_2(\tau)$
is uniformly bounded.
Suppose that this in not the case and this difference tends to infinity for some convergent sequence of $\tau_n \in [0,\tau_*]$:  $\tau_n \to \tau_0$.
Set
$$ v_n(x) = w_{\tau_n}(x+x_2(\tau_n)) . $$
Then $v_n(0)=w_2$, $v_n(x_1(\tau_n)-x_2(\tau_n))=w_1$.
We can choose a locally convergent subsequence from the sequence $v_n(x)$. Denote its limit by $v_0(x)$.
Then it is a solution of equation (\ref{b1}), $v_0(0)=w_2$, $v_0(-\infty) \leq w_1$ since
$x_1(\tau_n)-x_2(\tau_n) \to -\infty$. Then $v_0(-\infty)=w_0$ and therefore $c<0$ in view of Lemma 4.3.

Similarly, consider the sequence
$$ z_n(x) = w_{\tau_n}(x+x_1(\tau_n)) . $$
Then $z_n(0)=w_1$, $z_n(x_2(\tau_n)- x_1(\tau_n))=w_2$.
We can choose a locally convergent subsequence from the sequence $z_n(x)$. Let $z_0(x)$ denote its limit.
Then it is a solution of equation (\ref{b1}), $z_0(0)=w_1$, $z_0(+\infty) \geq w_2$ since
$x_2(\tau_n)-x_1(\tau_n) \to \infty$. Then $z_0(-\infty)=1, \ z_0(+\infty)=w_0$ and by virtue of Lemma 4.3, $c>0$.

 The obtained contradiction ($0<c<0$) shows that functions $x_1(\tau)$ and $x_2(\tau)$
are uniformly bounded.  Hence, for some fixed $T$ independent of $\tau$, it holds
\begin{equation}\label{ibw}
w_\tau(x) \leq w_2, \ t \geq T-c_*\tau_*, \quad w_\tau(x) \geq w_1, \ t \leq - T+c_*\tau_*, \ \mbox{for all } \ \tau \in [0, \tau_*].
\end{equation}
A useful consequence of this result is compactness of the set $\mathfrak{T}$ of all `admissible' delays and speeds:
$$\mathfrak{T} = \{(\tau,c) \in [0,\tau_*] \times[-c_*,c_*]: \mbox{system}\ (\ref{b1}), (\ref{b4})\ \mbox{has a monotone wave  for these} \ \tau, c \}.$$
Indeed, if $w_{\tau_n}(x), \ w_{\tau_n}(0) \in [\alpha, \beta]$,  is a sequence of monotone bistable waves propagating with speeds $c_n$, then each converging subsequence of $w_{\tau_n}(x), \tau_n, c_n,$ has a limit also satisfying, in view of  inequalities (\ref{ibw}), the boundary conditions  (\ref{b4}).

We claim that  if the parameters  $1-w_1$ and $w_2$ are sufficiently small then the functions $u_\tau = w_\tau - \psi, \ w_\tau'$ admit exponential estimates
for $x \leq -T$ and $x \geq T$ which are uniform uniform with respect to $\tau \in [0,\tau_*]$.  This means that  there exist some positive numbers $K_1, \gamma_1$ independent on $\tau$ such that
\begin{equation}\label{exe}
 |1-w_\tau(x)| + |w'_\tau(x)| \leq K_1e^{\gamma_1 x}, \ x \leq -T; \quad  |w_\tau(x)| + |w'_\tau(x)| \leq K_1e^{-\gamma_1 x}, \ x \geq T.
\end{equation}
These estimates follow from the assumptions $f(0) >0, \ f'(1) \in (-1,0]$, and their proof is based on arguments adapted from  the exponential dichotomy theory \cite{cop,H}. Since the set of all admissible delays $\mathfrak{T}$ is compact,  it suffices to establish (\ref{exe}) locally, i.e. to prove that  estimate (\ref{exe}) is true within  sufficiently small neighbourhood ${\mathcal O}$ of each point $(h_0, c_0) \in [0,c_*\tau_*] \times [-c_*,c_*]$ (so that $K_1$ and $\gamma_1$ may depend on $(h_0,c_0)$).
In the appendix, we briefly outline the proof of the first inequality in (\ref{exe})  in   more  difficult case when  $c_{\tau_0} \geq 0$ and  $c_\tau>0$.     Finally, since functions $w_\tau(x), w'_\tau(x)$ are uniformly (in $x$ and $\tau$) bounded,  Lemma 4.4 is an immediate consequence of inequalities (\ref{exe}).
\hfill  $\Box$

\vspace{3mm}

\noindent {\bf Corollary 4.5.}
{\em Assume all the conditions of Lemma 4.4.  Then
the functions $u_\tau=w_\tau-\psi$ are uniformly (in $\tau \in [0,\tau_*]$) bounded in the norm of $C_\mu^{2+\alpha}(\mathbb R)$.}

\medskip

\noindent
{\bf Proof.}
We have
$$ (u+\psi)'' + c (u +\psi)' + w (1 - w - f(w(x+c\tau))) = 0 . $$
Set $v = u \mu$ and multiply the last equation by $\mu$. Then we obtain
\begin{equation}\label{s13}
  v'' + c v' + g(x, c, \tau) = 0 ,
\end{equation}
where
$$ g(x,c,\tau) = - 2 u'\mu'-u\mu'' - c u \mu' + \psi'' \mu + c \psi'\mu + w(1 - w - f(w(x+c\tau))) \mu . $$
In view of (\ref{exe}), this function is uniformly (in $\tau, c$) bounded in the norm of $C^1(\mathbb R)$.
%
%
%
%
Therefore the norm of $v$ in  $C^{2+\alpha}(\mathbb R)$ is also uniformly bounded due to
the Schauder estimate. \hfill  $\Box$

\medskip

\noindent {\bf Lemma 4.5.}
{\em Let $w_\tau:{\mathbb R}\to (0,1)$ be monotonically decreasing wave solutions of  equation $A_\tau(w_\tau - \psi) =0$
for possibly different values of $\tau \in [0,\tau_*]$.
 Then $u_\tau =w_\tau - \psi$ are uniformly bounded in the norm of $C_\mu^{2+\alpha}(\mathbb R)$.
}

\medskip

\noindent
{\bf Proof.} Clearly, $w_\tau(x)$ is a bistable wave of equation (\ref{m1}) propagating with the
speed $c(u_\tau)$.  In view of (\ref{s10}),  $|c| \leq c_*$ for some $c_* > 0$.

Suppose now that the set  $\{u_\tau: \tau \in [0,\tau_*]\}$ is not uniformly bounded in  $E_\mu= C_\mu^{2+\alpha}(\mathbb R)$. Then there exist sequences $\tau_n$ and $u_n: = u_{\tau_n}$ such that $c(u_n) \to c_\star$ and  $\|u_n\|_{E_\mu} \to +\infty$ as $n \to +\infty$. Moreover,  without loss of generality, we can assume that
$w^{(j)}_n(x): = w^{(j)}_{\tau_n}(x) \to w^{(j)}_\star(x), \ j =0,1,2,$ uniformly on compact subsets of $\mathbb R$.
Here $w_\star$ denotes some non-increasing bounded solution of (\ref{m1}) with $c = c_\star$.

We claim  that
 $w_\star(-\infty)=1, \   w_\star(+\infty)=0$.  Indeed, if $w_\star(+\infty) >0$ then, applying the Fatou's lemma, we get the following contradiction:
 $$
+\infty = c(w_\star - \psi) \leq \liminf c(u_n) = \lim c_n = c_\star < c_*.
 $$
On the other hand, if $w_\star(-\infty)=0$ then $w_\star\equiv 0$ so that $w'_n(0) \to 0$ and  $w_n(x) \to 0, \ n \to +\infty,$ uniformly on each half-line $[s, +\infty)$. In addition, since each $w_n(x)$ satisfies
\begin{equation}\label{m1n}
 w'' + c_nw' + a_n(x)w  = 0, \quad a_n(x):= 1-w_n(x) - f(w_n(x+c_n\tau_n)),
\end{equation}
where $a_n(x) \to 1- f(0)$ uniformly on $[0,+\infty)$, we can conclude  (e.g. see \cite[Proposition 1, p.34]{cop}) that, for some positive constant $K, \gamma$, it holds
$$
|w_n(x)| \leq Ke^{-\gamma t} (w_n(0)+ |w'_n(0)|), \  t \geq 0, \quad n =1,2,\dots
$$
This implies, however, that $c(u_n) \to -\infty$, in view of the Lebesgue's dominated convergence theorem.
The obtained contradiction shows that  $w_\star(-\infty) \in \{w_0,1\}.$  For a moment, let suppose that
$w_\star(-\infty)= w_0$.  In this case, for each positive $\epsilon$, the intervals $T_n(\epsilon): = \{x: w_0-\epsilon < w_n(x) < w_0 + \epsilon\}$ have lengths $d_n(\epsilon)$ converging to $+\infty$ as $n \to +\infty$. Consequently, the intervals $Q_n:=\{x: w_0-\epsilon < w_n(x) \leq (w_0 + 1)/2\}$ have lengths $q_n > d_n$.  If $x_n$ denotes the unique solution of equation $w_n(x_n) = (w_0 + 1)/2$, then the sequence of shifted waves $\{w_n(x+x_n)\}$ converges, uniformly on compact subsets of $\mathbb R$, to a bounded non-increasing solution $w_l: {\mathbb R} \to [0, 1]$ of  equation (\ref{m1}) considered with $c=c_\star$.  Since $w_l(0)= (w_0 + 1)/2$ and $w_l(x) \geq w_0-\epsilon$ for all $x \in [0, q_n]$ with $q_n \to +\infty$, we conclude that $w_l(-\infty)=1$ and $w_l(+\infty)=w_0$.
However, the simultaneous existence of non-increasing waves $w_l$ and $w_\star$ of  equation (\ref{m1}) considered with the same speed $c=c_\star$ contradicts conclusions of Lemma 4.3.

Consequently,   $w_\star(-\infty)=1, \   w_\star(+\infty)=0$ so that $w_n(0) \to w_\star(0) \in (0,1), \ n \to +\infty$. This means that $w_n(0)\in [\alpha, \beta] \in (0,1), \ n =1,2,3\dots $ for some appropriate fixed $\alpha, \beta$.  By Corollary 4.5,  $u_n, \ n=1,2,3\dots$ must be uniformly bounded in the norm $C_\mu^{2+\alpha}(\mathbb R)$ that contradicts our choice of this sequence. \hfill $\square$


\section{Existence of solutions}

Once topological degree is defined and a priori estimates of solutions are obtained,
we can use the Leray-Schauder method.

\medskip

\noindent
{\bf Theorem 5.1.}
{\em Suppose that $f\in C^4({\mathbb R}_+)$ satisfies conditions (\ref{m4}), (\ref{m5}), (\ref{m6}) and, in addition,  $f'(w) < 0$ for $w_* \leq w < 1$, where $w_*= f^{-1}(1)$.
Then problem (\ref{b1}), (\ref{b4}) has a monotonically decreasing solution for any
$\tau \geq 0$.}

\medskip

\noindent
{\bf Proof.}   Fix an arbitrary positive number $\tau_*$.
With operator $A_\tau: E_\mu \to F_\mu$ defined  in Section 2, we consider all solutions of the equation
\begin{equation}\label{s14}
  A_\tau(u) = 0, \quad \tau \in [0,\tau_*],
\end{equation}
such that the functions $w = u + \psi$ are monotonically decreasing.
We denote such solutions by $u_M$.
It follows from Lemma 4.5 that $ \| u_M \|_{E_\mu} \leq K , $
where a positive constant $K$ does not depend on solution and on $\tau \in [0,\tau_*]$.
If the function $w = u +\psi$ is not monotonically decreasing, then
we denote such solutions $u_N$.

We claim that there exists a positive constant $r$ such that
\begin{equation}\label{s15}
  \|u_M - u_N\|_{E_\mu} \geq r
\end{equation}
for any solutions $u_M$ and $u_N$.
Indeed, if this is not the case, then there are two sequences of solutions $u_M^i$ and $u_N^i$ such that
\begin{equation}\label{s16}
  \|u_M^i - u_N^i\|_{E_\mu} \to 0 ,  \;\;\; i \to \infty.
\end{equation}
 Therefore, since  $|u_M^i|_{E_\mu} \leq K$ for all $i$ and
 $A_\tau(u)$ is a proper operator with respect to $(u,\tau)$ (see Theorem 2.4), we can choose convergent subsequence
of solutions. Without loss of generality we assume that $u_M^i \to u_0$ in $E_\mu$, $\tau_i \to \tau_0$.
Then the function $w_0=u_0+\psi$ is a solution of problem (\ref{b1}), (\ref{b4}) for some $\tau_0$ and $c$.
Hence $w_0'(x) \leq 0$, and by virtue of Lemmas 3.1, 3.2, we have $w'(x) < 0$ for all $x \in \mathbb R$.

Set $w_i = u_N^i + \psi$. It follows from (\ref{s16}) that $\|w_i - w_0\|_{C^2} \to 0$ as $i \to \infty$.
This convergence contradicts Lemma 3.4 and proves (\ref{s15}).

We construct an open  ball of the radius $r/2$ around each solution $u_M$. Since the set of solutions is compact,
we can choose a finite subcovering of the set of solutions by the balls. Denote this domain by $D$.
It contains all solutions $u_M$ and it does not contain solutions $u_N$.

Consider the topological degree $\gamma(A_\tau,D)$. Since $A_\tau(u) \neq 0$ on $\partial D$, then
this degree does not depend on $\tau$. It remains to verify that it is different from $0$.
Indeed, in the non-delayed case ($\tau=0$), equation (\ref{s14}) has a unique solution $u_0=w_0-\psi$, where $w_0$ is the unique
solution of problem (\ref{b1}), (\ref{b4}) with $\tau=0$. The value of the degree $\gamma(A_0,D)$ equals
the index of this solution which can be found through the eigenvalues of the linearized operator.  
This computation  was already done in \cite[Chapter 3, \S 3.2 ]{VVV} where it was shown that $\gamma (A_0,D) =1$.  Hence, $\gamma(A_\tau,D) = 1$ for any $\tau \in [0,\tau_*]$. This proves the existence of solutions $u_M$ of  equation (\ref{s14}) for each $\tau \in [0,\tau_*]$.
\hfill   $\Box$

\medskip

We note that if we consider the case where $c>0$ or $c<0$, the conditions on the function $f(w)$ can be somewhat
weakened (see Lemmas 3.1, 3.2).
\section*{Appendix}
In this section, we briefly outline the proof of the first inequality in (\ref{exe})  in  the case when  $(h,c):= (\tau c_{\tau}, c),\ c_\tau>0$, belongs to some small neighbourhood  ${\mathcal O}$ of the point $(h_0,c_{\tau_0})$ with $c_{\tau_0} \geq 0$.  In such a case, it is convenient to transform equation (\ref{b1}) into usual delayed differential equation by inverting the time: $w(t) = v(- t)$.
In this way, instead of system (\ref{b1}), (\ref{b4}), we obtain
 \begin{equation}\label{s2-}
  v''(x) - c v'(x) + v(x) (1-v(x)-f(v(x-h))) = 0 , \;\;\; v(-\infty)=0, \;\; v(+\infty)=1.
\end{equation}
Thus, for $x \geq T$, the function  $u(t)= v(t)-1$ satisfies the following delay differential equation:
\begin{equation}\label{s2+}
  u''(x) - c u'(x) - (1+a(x))u(x)  -(f'(1)+ b(x))u(x-h) = 0,
\end{equation}
where, uniformly with respect to $x \geq T$, it holds
$$a(x)= u(x) = O(|1-w_1|), \ b(x) = -f'(1)+(1+u(x))f'(1+\theta u(x-h)) = O(|1-w_1|), \ \theta \in (0,1).$$
We will extend $a(x), b(x)$ continuously on the whole $\mathbb R$ in such a way that $|a|_\infty, |b|_\infty = O(|1-w_1|)$. In the standard way, equation (\ref{s2+}) generates a semi-flow on the extended phase space $\mathbb R\times {\mathcal X}$, where  ${\mathcal X}:=C\times \mathbb R$  and  $C$ stands for the space of scalar continuous functions $C[-c_*\tau_*,0]$ provided with the sup-norm $|\phi|_\infty = \sup\{\phi(s): s \in [-c_*\tau_*,0]\}$.
In view of the assumption $f'(1)\in (-1,0]$, the `$\omega$-limit' equation of (\ref{s2+}),
\begin{equation}\label{cc}
u''(x) - c_0 u'(x) - u(x)  -f'(1)u(x-h_0) =0,
\end{equation}
has only the trivial bounded solution and therefore it possesses an exponential dichotomy with some projection $P_{h_0}(s): {\mathcal X} \to {\mathcal X}, \ s \in \mathbb R,$ and positive constants $K_0, \gamma_0$. The latter amounts to the following two properties:
\begin{itemize}
\item   If $(\phi_+, a_+) = P_{h_0}(s)(\psi,b)$,  with $(\psi, b)$ being an arbitrary fixed element of  ${\mathcal X}$,  then the solution $u(t, s,\phi_+,a_+),\ t \geq s,$ of the initial value problem $u(s+x)= \phi_+(x),$ $x \in [-c_*\tau_*,0]$,
$u'(s) = a_+,$ for equation (\ref{cc}) satisfies the inequality
$$
|u(t+\cdot, s, \phi_+,a_+)|_\infty+|u'(t, s, \phi_+,a_+)| \leq K_0e^{-\gamma_0(t-s)}(|b| + |\psi|_\infty), \quad t \geq s.
$$
\item On the other hand, if $(\phi_-, a_-) = (\psi,b) - P_{h_0}(s)(\psi,b)$ then the solution $u(t, s, \phi_-,a_-)$ of the initial value problem $u(s+x)= \phi_-(x), \ x \in [-c_*\tau_*,0]$,
$u'(s) = a_-,$ for equation (\ref{cc}) can be extended for all $t \leq s$ and satisfies the inequality
$$
|u(t+\cdot, s, \phi_-,a_-)|_\infty+|u'(t, s, \phi_-,a_-)| \leq K_0e^{\gamma_0(t-s)}(|b| + |\psi|_\infty), \quad t \leq s.
$$
\end{itemize}
 Observe that since equation (\ref{cc}) has constant coefficients,  $P_{h_0}(s)$ is also a constant function,  $P_{h_0}(s)\equiv P_0$.  We claim that $P_0(0,1) \not= (0,1)$.  Indeed, if $P_0(0,1) = (0,1)$, then the solution $\hat u(t):= u(t, 0,0,1),\ t \geq 0,$ of equation (\ref{cc}) is exponentially converging to $0$ as $t\to +\infty$ while $\hat u(s) =0$ for
$s \in [-h_0,0]$ and  $\hat u'(0) =1$.   This means that $\hat u(t)$  reaches its positive absolute maximum at some
leftmost point $t_M >0$. At this point,  $\hat  u''(t_M) \leq 0$, $\hat  u'(t_M) =0$,  so that, taking into account the inequality $0 \leq -f'(1) < 1$,  we get the following contradiction
$$
\hat  u(t_M) = \hat  u''(t_M) - c_0 \hat  u'(t_M)  -f'(1)\hat  u(t_M-h_0) \leq -f'(1)\hat  u(t_M-h_0) < \hat  u(t_M).
$$

Next, the roughness property of the exponential dichotomy guarantees (cf. \cite[Theorem 7.6.10]{H}) the existence of small $\delta >0$ such that
equation (\ref{s2+}) possesses an exponential dichotomy with some projection $P_{h}(s): {\mathcal X} \to {\mathcal X}, \ s \in \mathbb R,$ and  constants $2K_0, 0.5\gamma_0$ for all non-negative $c,h, w_1$ such that
\begin{equation}\label{nei}
\max\{|c-c_0|, |h-h_0|, |w_1-1|\} < \delta.
\end{equation}
Moreover, $\sup_{s \in \mathbb R}|P_0 - P_{h}(s)| \to 0$ as $\delta \to 0$. In particular, since $(Id -P_0)(0,1)\not = (0,0)$, we can take $\delta$ sufficiently small to have $i_* = \inf_{s \in \mathbb R}|(Id - P_{h}(s))(0,1)| >0$ once
(\ref{nei}) is satisfied.

Hence, assuming (\ref{nei}) and taking arbitrary solution $u(t) < 1-w_1, \ t \geq T-h,$ $u(+\infty)=$ $ u'(+\infty)=0$  of equation (\ref{s2+}),  we can conclude that $P_h(t)(u(t+\cdot), u'(t)) = (u(t+\cdot), u'(t))$ so that, for all $t \geq T$,
$$
i_* |u'(t)| \leq   |u'(t)| |(P_h(t) - Id)(0, 1)|= |(Id- P_h(t))(u(t+\cdot), 0)| \leq  |P_h(t) - Id| |u(t+\cdot)|_\infty.
$$
Thus we obtain that, for all $c,h$ satisfying (\ref{nei}) there exists some universal constant $C>0$ such that
$
|u'(t)| \leq  C  |u(t+\cdot)|_\infty, \ t \geq T.
$
This yields the required uniform exponential estimate
$$
|u(t+\cdot )|_\infty+|u'(t)| \leq 2K_0e^{-0.5 \gamma_0(t-T)}(|u(T+\cdot )|_\infty+|u'(T)|) \leq $$
$$ 2(1+C)K_0e^{-0.5 \gamma_0(t-T)}|u(T+\cdot )|_\infty \leq 2(1+C)(1-w_1)K_0e^{-0.5 \gamma_0(t-T)}, \quad t \geq T,
$$
which holds for all $c\geq 0,h$ satisfying (\ref{nei}).
\section*{Acknowledgements}

The first author was supported by FONDECYT (Chile), project 1150480.
The second author was partially supported  by the Ministry of Education and Science of the Russian Federation (the Agreement number 02.a03.21.0008).


\end{document}